\documentclass[11pt]{amsart}

\usepackage[letterpaper,hmargin=1.4in,vmargin=1.2in]{geometry}
\usepackage{amsmath, amssymb, amsthm, verbatim, amscd, graphicx}

\DeclareSymbolFont{AMSa}{U}{msa}{m}{n}
\DeclareMathSymbol{\square} {\mathord}{AMSa}{"03}

\DeclareSymbolFont{AMSb}{U}{msb}{m}{n}
\DeclareMathSymbol{\N}{\mathbin}{AMSb}{"4E}
\DeclareMathSymbol{\Z}{\mathbin}{AMSb}{"5A}
\DeclareMathSymbol{\R}{\mathbin}{AMSb}{"52}
\DeclareMathSymbol{\Q}{\mathbin}{AMSb}{"51}
\DeclareMathSymbol{\I}{\mathbin}{AMSb}{"49}
\DeclareMathSymbol{\C}{\mathbin}{AMSb}{"43}

 \newcommand{\ga}{{\mathfrak g}}
  \newcommand{\gh}{\hat{\mathfrak g}}
  \newcommand{\ghh}{\widehat{\mathfrak{sl}}(2) }

\begin{document}

\title{Families of Dirac operators and quantum affine groups}

\author{Jouko Mickelsson}

\address{Department of Mathematics and Statistics,
  University of Helsinki.
 Department of Theoretical Physics,
  Royal Institute of Technology, Stockholm}

\date{June 21, 2010,  "Geometry and Quantum Field Theory" / Careyfest  60/ MPIM, Bonn}

\maketitle

\bf Abstract \rm Twisted K-theory classes over  compact Lie groups can be realized as families of Fredholm operators  using the representation theory
of loop groups.  In this talk I want to show how to deform the Fredholm family, in the sense of quantum groups.  The family of Dirac type operators is 
parametrized by vectors in the adjoint module for a quantum affine algebra and transform covariantly under a (central extension of) the algebra.

\section{Introduction}

Let $X$ be  a topological  space, $Fred_*$ the space of self-adjoint Fredholm
operators in a complex Hilbert space $H$ with both positive and negative essential spectrum. 
This is a universal 
classifying space for $K^1$. Actually, one can take as the \bf definition: \rm

$$K^1(X) = \{\text{homotopy classes  of maps }f:X \to Fred_* \}$$
In the even case
$$K^0(X) = \{\text{homotopy classes of maps } f:X\to Fred \}$$
where $Fred$ is the space of all Fredholm operators in $H.$

The Chern character 
$$ch: K^1(X) \to H^{odd}(X, \Bbb Z)$$
is an additive map to odd cohomology classes.  In particular, the degree 3 component
$DD(f)=ch_3(f)$ of $[f]  \in  K^1(X)$ is called the \bf Dixmier-Douady class \rm of the gerbe defined 
by the the family $f(x)$ of Fredholm operators. 
In the de Rham cohomology an equivalent construction of $DD(f)$ comes from the family $L_{\lambda\lambda'}$ of complex
line bundles. One can choose the curvature forms $\omega_{\lambda\lambda'}$ such that
$$\omega_{\lambda\lambda'} + \omega_{\lambda'\lambda''} = \omega_{\lambda\lambda''}$$
and with a partition of unity $\sum \rho_{\lambda}=1$ subordinate to the cover by the open sets
$U_{\lambda}$ one has
$$DD(f) = \sum_{\lambda}  d\rho_{\lambda}\wedge \omega_{\lambda\lambda'} $$
and this does not depend on the choice of $\lambda'.$ The 3-cohomology class is related to quantum field theory anomalies 
and it transgresses in the Hamiltonian formulation of gauge theory to central, or more generally, abelian extension
of the group of gauge transformations, [CMM].

A 3-cohomology class $[H]$  can be used as twisting in K-theory. Instead of a globally defined family of Fredholm operators
over $X$ we have local families $f_{\alpha}: U_{\alpha} \to Fred(H),$ where $\{U_{\alpha}\}$ is an open cover of $X,$
such that on the overlaps $U_{\alpha} \cap U_{\beta}$ one has
$$ f_{\beta}(x) = h_{\beta\alpha}(x) f_{\alpha}(x) h_{\beta\alpha}(x)^{-1}$$
where the $h_{\alpha\beta}$'s are the transition functions for a principal $PU(H)$ bundle $P$ over $X.$ 
The equivalence class of the bundle $P$ is given by $[H].$ The local 2-forms $\omega_{\alpha\beta}$ are given as
the pull-backs $h_{\alpha\beta}^*\omega$ where $\omega$ is the Chern class of the canonical complex line bundle over
$PU(H).$ 

Quantum field theory provides a method for constructing twisted families of Fredholm operators. In particular, one can use the
supersymmetric Wess-Zumino-Witten model in the case when $X$ is a compact Lie group. This construction actually provides 
an operator theoretic realization for the relation between twisted K-theory and the Verlinde algebra in conformal field theory,
proven in [FHT]. In the next section I will briefly recall the construction of the undeformed family of Dirac operators, [M].

\section{The supersymmetric WZW model}
Families of Dirac operators $D_A$ on the unit circle $S^1$ coupled to smooth vector potentials $A$ transform covariantly under the  gauge group
action, defining an element in $K^*(\mathcal A/\mathcal G)= K^*(G).$  Here $\mathcal G =\Omega G,$ the based  smooth loop group of a compact Lie group $G.$ 
The quantized Dirac
operators $\hat D_A$ acting in a fermionic Fock space transform covariantly under a central extension $\widehat{LG}$ of the full loop group $LG,$ thus the family
$\hat D_A$ is a candidate for an element in twisted K-theory on the moduli space of gauge connections. However, there is a catch: These operators
are essentially positive, we need operators with both positive and negative essential spectrum.
 
To find  genuine nontrivial twisted K-thery classes we have to go to the supersymmetric WZW model. Morally, the family of Fredholm operators is
now a family of Dirac operators on the loop group $LG.$ We cannot make sense of the Dirac operators on the infinite-dimensional manifold $LG$ 
by the standard analytic methods (we would need a Haar measure on the group $LG$) but there is a purely algebraic construction using the highest weight
representation theory of $\widehat{LG}.$ 

Concretely, the operators are acting in a tensor product $H= H_f \otimes H_b,$ where $H_b$ carries an irreducible highest weight representation of
$\widehat{LG}$ of level $k=0,1,2,\dots$ and $H_f$ carries an irreducible representation of the Clifford algebra of the vector space $L\mathfrak g.$ The 
Clifford algebra is then used to construct a (reducible) representation of the central extension $\hat{\mathfrak g}$ of $L\mathfrak g,$ of level $\kappa,$ the
dual Coxeter number of $\mathfrak g.$ A basis for the Clifford algebra consists of elements $\psi^n_a$ with $n\in\mathbb Z$ and $a=1,2,\dots, \text{dim}\, \mathfrak g$
with anticommutation relations 
$$\psi_a^n \psi_b^m + \psi_b^m \psi_a^n = 2 \delta_{ab} \delta_{n,-m}.$$
The commutation relations of the basis vectors $T^n_a$ of $L\mathfrak g$ in the representation space $H_b$ are
 $$ \left[T_a^n, T_b^m \right] = \lambda_{abc} T_c^{n+m} + \frac{k}{4} \delta_{ab} n\delta_{n,-m}$$
where the $\lambda_{abc}$'s are the structure constants of $\mathfrak g$ in a suitably normalized basis
 and the Dirac operator $Q_A$ coupled to a vector potential on $S^1$ is
\begin{eqnarray}
Q_A &= &i\psi_a^n \otimes T_a^{-n} - \frac{i}{12} \lambda_{abc} \psi_a^n \psi_b^m\psi_c^{-n-m} \otimes 1
+ i\frac{k+\kappa}{4} \psi_a^n A_a^{-n} \otimes 1 \\
&=& i\psi_a^n \otimes T_a^{-n} +  \frac{i}{3} \psi^n_a K_a^{-n}\otimes 1  + i\frac{k+\kappa}{4} \psi_a^n A_a^{-n}\otimes 1,
\end{eqnarray}
where $K_a^n$'s satisfy the same commutation relations, as operators in $H_f,$ as the operators $T^n_a$ except that the level is $\kappa$ instead of $k.$
Here $A_a^n$'s are the Fourier components of a vector potential on the circle. 

The family $Q_A$ 
transforms covariantly under the projective representation of level $k+\kappa$ the loop group $LG$ defining an element in $K(G, k+\kappa)$ corresponding to the D-D class $[H]$ in $H^3(G,\mathbb Z)$
equal to $k+\kappa$ times the basic class in $H^3(G) = \mathbb Z$ when $G$ is a simple simply
connected compact Lie group.

Actually, since $\mathcal A/ \Omega G = G$ and $G\subset LG$, we have an $G$ equivariant class, element of
$K^*_G(G,  H).$  In the q-deformed case the gauge symmetry is not a group but a Hopf algebra, and the moduli space of 
gauge connections is not defined. Instead, the role of the Dixmier-Douady class $[H]$ is completely taken over by the (level of) 
the central extension of the quantum affine algebra.

Our construction can be viewed as a generalization of [BK], from the quantum $SU(2)$ algebra to quantum affine algebras. 
A more detailed article is in preparation, [HM].  I want to stress that we are not trying to construct a spectral triple, in the sense of
noncommutative geometry, using the generalized Dirac operator. Such constructions for compact quantum groups have been discussed in several papers during
the last ten years. I just refer to one recent paper [NT], see for references to earlier literature therein.

\section{Quantum affine algebra}

Let $\mathfrak g$ be a simple finite-dimensional Lie algebra and $\hat{\mathfrak g}$ the associated affine Lie algebra. The quantum affine algebra
$U_q(\gh)$ is generated by  vectors $e_0,e_1, \dots ,e_{\ell}, f_0,f_1,\dots,  f_{\ell}, \newline
K_0, K_1, \dots, K_{\ell}, K_0^{-1}, \dots, K_{\ell}^{-1}$ with
the relations

\begin{eqnarray*}
[e_i,f_i] &= & \delta_{ij} \frac{K_i -K_i^{-1}}{q - q^{-1}}, K_iK_j = K_j K_i  \\
K_i e_j K_i^{-1} &= & q^{\alpha_{ij}} e_j,  K_i f_j K_i^{-1} = q^{-\alpha_{ij}}f_j \\
\sum_{k=0}^{1-a_{ij}} &(-1)^k &\left[ \begin{array}{c} 1-a_{ij} \\ k \end{array} \right]_q e_i^{1-a_{ij} -k} e_j e_i^k = 0 \hskip 0.2in (i\neq j) \\
\sum_{k=0}^{1-a_{ij}}  &(-1)^k & \left[ \begin{array}{c} 1-a_{ij} \\ k \end{array} \right]_q f_i^{1-a_{ij} -k} f_j f_i^k = 0 \hskip 0.2in (i\neq j) 
\end{eqnarray*}

where 

$$ \left[\begin{array}{c} m\\k \end{array} \right]_q = \frac{ m_q (m-1)_q \dots (m-k+1)_q }{ k_q (k-1)_q \dots 1_q}$$
$$ k_q = 1 + q +\dots q^{k-1}.$$
Here the parameter $q$ is a positive real number and the integers $a_{ij}$ are the matrix elements of the Cartan matrix of $\gh.$

Let $A^n_i$ with $n\in \Bbb Z$ and $i$ is a label for the basis in the adjoint module of $\mathfrak g,$ be a basis for the q-affine adjoint module. Under $\ga$ each 'Fourier mode' $A^n$ transforms
according to the adjoint representation of $U_q(\ga)$, which is a q-deformation of the adjoint representation of $\ga.$ The generator $e_0$ increases the index $n$ by
one unit, $f_0$ decreases it by one unit. For example, for $\ga = \mathfrak{sl}(2)$ one has the explicit formulas
 
\begin{eqnarray*}
  e_1 A^n_1 &=& f_0 A^n_1=0, \hskip 0.4in  f_1 A^n_1  = A^n_0, \hskip 0.3in e_0 A^n_1 = A^{n+1}_0 \\ 
e_1 A^n_0  &=& (q+q^{-1}) A^n_1, \hskip 0.3in  f_0 A^n_0 = (q+q^{-1})A^{n-1}_1\\
 f_1 A^n_0 &= &A^n_{-1}, \hskip 0.7in  e_0 A^n_0 = A^{n+1}_{-1}\\
e_1 A^n_{-1} &=& (q+q^{-1}) A^n_0, \hskip 0.3in  f_0 A^n_{-1}=  (q+q^{-1}) A^{n-1}_0\\ 
f_1 A^n_{-1}&=&0= e_0 A^n_{-1}\\
K_1 A^n_i &= &q^{2i} A^n_i = K_0^{-1} A^n_i. \end{eqnarray*}

The vectors $A^n_i$ will be constructed as operators acting in a Fock space carrying a representation of $U_q(\gh)$ 
such that the adjoint action is given by
$$x.A^n_i  = \sum_{(x)} x'A^n_i S(x'')   \text{ for } x\in U_q(\gh),$$
where $S:U_q(\gh) \to U_q(\gh)$ is the antipode and $\Delta(x) = \sum_{(x)} x' \otimes x''$ is the coproduct $\Delta: U_q \to U_q \otimes U_q.$ 
We also need the Clifford algebra generated by elements $\psi^n_i$ acting in the Fock space and transforming under $U_q(\gh)$ according to the dual
adjoint representation (which in fact is equivalent to the adjoint representation).

\section{The Dirac operator}

The Dirac operator $Q$ is acting in $H_f \otimes H_b$ where $H_f$ is 
the q-fermionic Fock space and $H_b$ carries another highest weight representation of $U_q(\gh).$  The action of the nontrivial central
extension is seen in the action of the elmenent $K_0K_1\dots K_{\ell},$ which is no more equal to the unit operator but a power of $q$ where
the exponent depends on the level of the representation.
\begin{eqnarray*}
Q = i \sum \psi^n_a \otimes T^{-n}_a + i \frac13 \sum  \psi^n_a K^{-n}_a \otimes 1
\end{eqnarray*}
where $T^n_a$ are basis vectors of the adjoint module, acting as linear operators in the space $H_b.$  We need also another copy of the adjoint module,
acting in the space $H_f.$ The components are denoted by $K_a^n.$ The vectors $\psi^n_a$ are elements in a quantum Clifford
algebra acting as operators in a q-Fock space $H_f,$ to be defined below in Subsection 4.1.

In contrast to the undeformed case, the operators $K_a^n, T_a^n$ do not satisfy the defining relations of the algebra $U_q(\gh).$

Let $R$ be the universal R-matrix for the algebra $U_q.$ An explicit construction is given in [KT]. 
Following [DG], we can then define a basis for vectors in a submodule $A\subset  U_q$ transforming according 
to an adjoint representation 
$$ ad_q(x) v = \sum_{(x)} x' v S(x'')$$
of $U_q$ on itself.    A basis is defined as
$$ A^n_i = \sum K_{n,i}^{m,\alpha; p,\beta}(\pi_{m,\alpha;p,\beta} \otimes id) A,$$
where $A =(R^T R -1)/h,$ with $e^h =q$ and $R^T = \sigma R\sigma,$ where $\sigma$ permutes the factors in the
tensor product $U_q \otimes U_q.$  Here $\pi_{m,\alpha;p,\beta}$ are the matrix elements in the defining representation $V$
of $U_q.$

For example, for $\gh = \widehat{\mathfrak{sl}}(2)$ the basis in the defining representation is $v^n_i$ with $n\in\Bbb Z$ and $i=-1,0,1$ and $\alpha,\beta= \pm.$ The numerical coefficients
$K$ come from the identification of the basis of the adjoint representation as linear combinations of the basis vectors in $V\otimes V.$ 

The action of the 
Serre generators in the defining representation is
\begin{eqnarray*}
 e_1 v^n_+  &=& f_0 v^n_+ = 0, f_1 v^n_+ = v^n_-, e_0 v^n_+ = v^{n+1}_-, e_0 v^n_+ = v^{n-1}_-  \\
e_1 v_-^n &=& v^n_+,  \,\,  f_0 v^n_- = v^{n-1}_+, e_0 v^n_- =0= f_1 v^n_-  \\
K_1 v^n_{\pm} &= & q^{\pm 1} v^n_{\pm} = K_0^{-1} v^n_{\pm}.      \end{eqnarray*}

\subsection{Generalized affine Hecke algebra, the case of  $U_q(\ghh)$ }

The affine Hecke algebra for $\gh$ is defined through the relations coming from  the R-matrix
 $\text{\v R} = \sigma R$ in the tensor product $V^0 \otimes V^0.$  See [L] for the use of the Hecke algebra for constructing
 Fock space representations of $U_q(\widehat{sl}(n)).$ The matrix satisfies
$$ (\text{\v R} - q^{-1})(\text{\v R} +q) = 0,$$
since $-q$ and $q^{-1}$ are the only eigenvalues of the invertible matrix \v R. Denote by $Y_1$ the shift operator which sends $v^n_i \otimes v^m_j$
to $v^{n+1}_i \otimes v^m_j$ and by $Y_2$ the corresponding shift operator acting on the second tensor factor. 
The matrix \v R acting on $V$ is then defined using the relations 
$$ \text{\v R} Y_1 = Y_2 \text{\v R}^{-1},  \,\,   \text{\v R} Y_2 =  Y_1 \text{\v R} +(q-q^{-1}) Y_2.$$
Actually, the second relation follows from the first and the minimal polynomial relation.

Now the braiding relations are given by setting the ideal in the tensor algebra of $V$ generated by the elements
$$ (q^{-1} + \text{\v R})(V\otimes V) $$
equal to zero. These have in particular the consequence that any $v^n_i v^m_j$ with $n > m$ can be written as a linear combination of vectors
$v^p_k  v^q_l$ with $p+q = n+m$ and $p\leq q.$  In the zero mode space $V^0$ the meaning of the braiding relations is that they project out the 
'symmetric' part of the tensor product $V^0\otimes V^0.$ The 3-dimensional representation
is the eigenspace of \v R with eigenvalue $q^{-1}$ and the 1-dimensional component corresponds to the eigenvalue $-q.$

To complete the construction of the Dirac operator we need also the generalized Clifford algebra in the coadjoint representation. The algebra is generated 
by vectors $\psi^n_i$ with $n\in \Bbb Z$ and $i=1,0, -1.$ The defining relations are given by braiding relations and an invariant (nonsymmetric) bilinear 
form. The braiding relations are defined recursively like in the case of $V$ (or $V^*$)  with the difference that since the R-matrix $\text{\v R}$ in the adjoint representation
has 3 instead of 2 different eigenvalues, which are now $-q^{-2}, q^2, q^{-4},$ with multiplicities $3,5,1$ respectively. 
The Fock space representation of the Clifford algebra is defined similarly to the undeformed case: It is generated by a finite-dimensional irreducible spin representation
of the zero Fourier mode algebra, with the additional relation that  this subspace is annihilated by the elements $\psi^n_i$ for $n<0.$ 

 Finally, the action of $U_q(\widehat{\mathfrak{sl}}(n))$
is defined in the tensor algebra modulo the ideal generated by the elements (3) below, over the coadjoint module, using the opposite coproduct $\Delta^{op}(x)= \sum_{(x)} x'' \otimes x'.$

The negative eigenvalue corresponds again
to a 3-dimensional 'antisymmetric' representation and the positive eigenvalues to a 6-dimensional 'symmetric' representation; the latter contains the 
1-dimensional trivial representation. 

The Hecke algebra is replaced by a generalized Hecke algebra,
\begin{eqnarray*} 
 & & Y_1 Y_2 = Y_2 Y_1  \\ 
& &(\text{\v R} - q^2)(\text{\v R} -q^{-4})(\text{\v R} + q^{-2}) =0 \\
& & \text{\v R} Y_1 = Y_2 \text{\v R}^{-1} ,  \,\, \text{\v R} Y_2 = Y_i  \text{\v R} +(q^2 -q^{-2}) Y_2
\end{eqnarray*}
where the middle relation is the minimal polynomial of the diagonalizable matrix \v R.

The generalized symmetric tensors correspond to positive eigenvalues of  \v R. In the Clifford algebra symmetrized products are identified as scalars
times the unit.  That is, we fix a $U_q(\ghh)$ invariant bilinear form  $B$  and the Clifford algebra is defined as the tensor algebra over $V$ modulo the ideal 
generated by 
\begin{eqnarray}
P(u \otimes v) - B(u,v)\ \cdot 1 
\end{eqnarray}
where $P$ is the projection on positive spectral subspace of \v R.  In the case when $V$ is the adjoint module for $U_q(\ghh)$ one can fix $B$ by 
identifying the first factor $V$ as the dual $V^*$ and using the natural pairing $V^* \otimes  V \to \Bbb C.$ 
Alternatively, one can view $B$ as the projection onto the 1-dimensional trivial submodule inside of the 'symmetric module'.

\subsection{The family of Dirac operators}

In the nondeformed case one has for an infinitesimal gauge transformation $X\in L\mathfrak{g}$
$$ [X, Q] =  (k + \kappa)\sum (-n) \psi^n_i X^{-n}_i = \frac{k+\kappa}{4} <\psi, dX> $$
and for a family of operators $Q_A = Q + \frac{k+ \kappa}{4}\psi^n_i A^{-n}_i$
$$[X, Q_A] =  \frac{k+\kappa}{4}<\psi, [A,X] +dX>.$$
In q-deformed case $A$ is to be understood as a vector in the adjoint module extended by $\Bbb C c.$ The above equation is replaced by
$$ \sum_{(x)} x' Q_A S(x'') =  Q_{x._c A}$$
where $x_c$ denotes the adjoint action in the centrally extended module,
$$x._c  A = x. A + \lambda_x(A) c$$
with $\lambda_x$ a linear form on the adjoint module, linear in the argument $x\in U_q(\ghh).$ 
It satisfies the cocycle relation
$$\lambda_{xy}(A) = \lambda_x(y.A).$$
Here $c$ is an element of the extended 
module such that $x.c = 0$ for all $x\in U_q(\ghh).$ It is enough to give value of $\lambda_x$ when $x$ is a Serre generator. The only nonzero
forms are $\lambda_{e_0}$ and $\lambda_{f_0}.$ The former is nonzero only for the component $A = A^{-1}_1$ and the latter for $A=A^1_{-1}.$
In the particular case when the adjoint module is acting as operators in a level = 1 representation of $U_q(\ghh)$ we have 
 $\lambda_{e_0} (A^{-1}_1) =
-q^{-1}$  and  $\lambda_{f_0}(A^1_{-1}) = q^{-1}.$ 

 \bf References\rm

[CMM] A.L. Carey, J. Mickelsson, and M. Murray:  Bundle gerbes applied to quantum field theory.
Rev. Math. Phys \bf 12, \rm 65 (2000)

[DG] Gustav W. Delius, Mark D. Gould, Andreas H\"uffmann, Yao-Zhong Zhang:
Quantum Lie algebras associated to $U_q(gl_n)$ and $U_q(sl_n).$
 arXiv:q-alg/9508013. \newline
 J.Phys.A29:5611-5618,1996

[FHT] Daniel S. Freed, Michael J. Hopkins, Constantin Teleman:
Twisted equivariant K-theory with complex coefficients.
arXiv:math/0206257. J. Topology; January 2008: 1, 16-44

[BK]  P.N. Bibikov, P.P. Kulish:  Dirac operators on quantum $SU(2)$ group and quantum sphere
arXiv:q-alg/9608012 

[KT]   S. Khoroshkin, V. Tolstoy. Twisting of quantum (super)algebras. Connection of DrinfeldÕs and Cartan-Weyl realizations for quantum affine algebras. 
 hep-th/9404036

[M] J. Mickelsson:  Gerbes, (twisted) K-theory, and the supersymmetric  WZW model. Strasbourg 2002,
 Infinite dimensional groups and manifolds,  93--107, IRMA Lect. Math. Theor. Phys., 5, de Gruyter, Berlin, 2004. 
 arXiv:hep-th/0206139
 
 [HM] Antti Harju and Jouko Mickelsson, in preparation
 
 [L] Bernard Leclerc: Fock space representations of $U_q(\widehat{sl(n)}).$  Lecture notes, Grenoble 2008,
 http://www-fourier.ujf-grenoble.fr/IMG/pdf/leclerc{\_}rev.pdf

[NT]  Sergey Neshveyev and  Lars Tuset: The Dirac Operator on Compact Quantum Groups. 
arXiv:math/0703161

\enddocument